\documentclass[article,11pt]{amsart}
\usepackage[top=23.86mm,bottom=24mm,left=25mm,right=25mm]{geometry}

\newtheorem{lemma}{Lemma}[section]
\newtheorem{proposition}{Proposition}[section]
\newtheorem{theorem}{Theorem}[section]
\newtheorem{corollary}{Corollary}[section]
\newtheorem{definition}{Definition}[section]
\newtheorem{example}{Example}[section]
\newtheorem{remark}{Remark}[section]

\makeatletter
\@addtoreset{equation}{section}
\makeatother

\DeclareMathOperator{\Spec}{Spec}

\DeclareMathOperator{\Id}{Id}

\DeclareMathOperator{\loc}{loc}

\providecommand{\pro}[1]{(#1_t)_{t \geq 0}}

\providecommand{\semi}[1]{\{#1_t: t \geq 0\}}

\newcommand{\cK}{\mathcal{K}}

\newcommand{\R}{\mathbf{R}}
\newcommand{\1}{\mathbf{1}}

\newcommand{\pr}{\mathbf{P}}

\newcommand{\ex}{\mathbf{E}}
\newcommand{\Rd}{\mathbf{R}^d}
\newcommand{\N}{\mathbf{N}}

\begin{document}
\title[Ground state domination properties for non-local Schr\"odinger operators]
{Contractivity and ground state domination properties for non-local Schr\"odinger operators}
\author{Kamil Kaleta, Mateusz Kwa\'snicki, J\'ozsef L\H orinczi}
\address{Faculty of Pure and Applied Mathematics \\ Wroc{\l}aw University of Technology
\\ Wyb. Wyspia{\'n}skiego 27, 50-370 Wroc{\l}aw, Poland}
\email{kamil.kaleta@pwr.edu.pl, mateusz.kwasnicki@pwr.edu.pl}

\address{Department of Mathematics \\ Loughborough University \\ Loughborough LE11 3TU, United Kingdom}
\email{J.Lorinczi@lboro.ac.uk}

\thanks{\emph{Key-words}: symmetric L\'evy process, non-local Schr\"odinger operator, Feynman-Kac semigroup, ground state
transformation, $L^p$-ground state domination, intrinsic hypercontractivity, intrinsic supercontractivity, (asymptotic) intrinsic
ultracontractivity \\
\noindent
2010 {\it MS Classification}: Primary 47D08, 60G51; Secondary 47D03, 47G20 \\
\noindent
Research supported by the National Science Center (Poland) grant on the basis of the decision
No. DEC-2011/03/D/ST1/00311 and by the Foundation for Polish Science.
}

\begin{abstract}
We study supercontractivity and hypercontractivity of Markov semigroups obtained via ground state transformation of non-local
Schr\"odinger operators based on generators of symmetric jump-paring L\'evy processes with Kato-class confining potentials. This
class of processes has the property that the intensity of single large jumps dominates the intensity of all multiple large jumps,
and the related operators include pseudo-differential operators of interest in mathematical physics. We refine these
contractivity properties by the concept of $L^p$-ground state domination and its asymptotic version, and derive sharp necessary and
sufficient conditions for their validity in terms of the behaviour of the L\'evy density and the potential at infinity. As a consequence,
we obtain for a large subclass of confining potentials that, on the one hand, supercontractivity and ultracontractivity, on the other
hand, hypercontractivity and asymptotic ultracontractivity of the transformed semigroup are equivalent properties. This is in stark
contrast to classical Schr\"odinger operators, for which all these properties are known to be different.
\end{abstract}

\maketitle

\baselineskip 0.5 cm

\bigskip

\section{Introduction}
In recent years Markov processes with jumps and non-local pseudo-differential operators, which are generators of these processes,
have received much attention in both pure and applied mathematics. They allow to model a variety of discontinuous random phenomena in
statistical mechanics, anomalous transport, laser optics and, via a Feynman-Kac-type representation, relativistic quantum theory
and quantum models with spin. The non-local character often poses intriguing problems, especially in areas where functional analysis,
PDE and probability theory meet.

In many interesting cases the specific models involve non-local Schr\"odinger operators based on generators of L\'evy processes
with jumps. Recent investigations include heat trace and spectral gap estimates \cite{AVB, bib:BYY, Kal2012}, gradient estimates of
harmonic functions \cite{K2013}, properties of radial solutions, ground states, eigenfunctions and eigenvalues
\cite{KKMS,K12,KKM,FLS2015,bib:KL15b,LM,bib:HL}, smoothing properties of evolution semigroups \cite{bib:KL15a, ChW2015, bib:KaKu},
properties of the associated transformed jump processes \cite{bib:KL, KL15c}, as well as applications in quantum theory
\cite{bib:LS, HL08, HIL, HLT, HHL, BSSch1}.

In this paper we focus on fundamental properties of semigroups $\{T_t: t\geq 0\}$ of non-local Schr\"odinger operators of the form
$H=-L+V$ by using functional integral techniques. Here $L$ is the generator of a symmetric jump L\'evy process in $\R^d$, $d \geq 1$,
with the property that all multiple large jumps are dominated under the L\'evy measure with density $\nu$ by a single large jump (which
we call a \emph{jump-paring L\'evy process}). This property proves to be a unifying concept including many jump processes and related
operators of interest \cite{bib:KL15a}. The term $V$ will be chosen to be a Kato-class confining potential (see Section 3.1 below for
precise definitions), thus the operator $H$ has a purely discrete spectrum, consisting of eigenvalues $\lambda_0 < \lambda_1 \leq \lambda_2
\leq ... \to \infty$ of finite multiplicities. The unique eigenfunction $\varphi_0$ corresponding to eigenvalue $\lambda_0$, called
\emph{ground state}, has a continuous, bounded and strictly positive version (in Lebesgue sense), which will be our choice throughout
below.

We define the \emph{ground state-transformed} (also known as \emph{intrinsic}) semigroup $\{\widetilde T_t: t\geq 0\}$ associated with
the non-local Schr\"odinger operator $H$ by
$$
\widetilde T_t f(x) = \frac{e^{\lambda_0 t}}{\varphi_0(x)} T_t (f \varphi_0)(x), \quad f \in L^2(\R^d,\mu), \ \ t \geq 0,
$$
where $\mu(dx) = \varphi_0^2(x)dx$. The semigroup $\{\widetilde T_t: t\geq 0\}$ is conservative and determines a right Markov process
with stationary distribution $\mu$, which we call \emph{ground state-transformed jump process} corresponding to $H$ \cite{KL15c, bib:KL}.

The aim of this paper is to study the contractivity properties of the operators $\widetilde T_t$ from $L^2(\R^d,\mu)$ to $L^p(\R^d,\mu)$ for $p \in
(2,\infty)$. Recall that the semigroup $\{\widetilde T_t: t\geq 0\}$ is called \emph{hypercontractive} (resp. \emph{supercontractive}) if $\widetilde
T_t$ maps $L^2(\R^d,\mu)$ into $L^p(\R^d,\mu)$ for every $p \in (2,\infty)$ and for some $t  >0$ (resp. for every $t>0$). Moreover, $\{\widetilde T_t:
t\geq 0\}$ is \emph{ultracontractive} (resp. \emph{asymptotically ultracontractive}) if the same holds with $p=\infty$ for every $t>0$ (resp. for some
$t>0$). Such smoothing properties have been widely studied for classical Schr\"odinger operators $H= -\Delta + V$ \cite{bib:S}. Hypercontractivity
has been introduced by Nelson \cite{N66}, and proved to be useful especially in quantum field theory. The first works on ultracontractivity are due
to Davies and Simon \cite{bib:DS}, and Ba\~nuelos \cite{bib:B}, and it has a number of useful consequences (for a detailed discussion see
\cite[Sect. 1]{bib:KL15a}). Generally, such contractivity properties are useful in establishing functional (such as log-Sobolev, Nash, Poincar\'e
etc) inequalities, as well as heat kernel estimates \cite{bib:BGL,bib:D}. For the classical case there is a natural hierarchy of these properties.
For instance, in the case of $H=-\Delta + |x|^{\alpha} (\log(1+|x|)^{\beta}$, the related transformed semigroup is not hypercontractive if $\alpha < 2$,
it is hypercontractive but not supercontractive if $\alpha = 2$, $\beta = 0$, it is supercontractive but not ultracontractive if $\alpha = 2$, $0 < \beta
\leq 2$, and it is ultracontractive if $\alpha = 2$, $\beta  > 2$, or if $\alpha > 2$. This shows that ultracontractivity is essentially stronger than
supercontractivity, which is essentially stronger than hypercontractivity. Note that the asymptotic version of ultracontractivity has not been studied
in the classical case before the paper \cite{bib:KL}, where it was proposed first.

For non-local Schr\"odinger operators the situation is less explored (see \cite{ChW2015}), with the exception of ultracontractivity properties, i.e.,
when $p=\infty$. This has been studied in the recent paper \cite{bib:KL15a} for the same class of underlying jump-paring processes and Kato-class
potentials (see this reference also for a discussion of previous literature). An interesting feature of non-local Schr\"odinger operators is that
all eigenfunctions are pointwise dominated everywhere by the ground state \cite[Cor. 2.1, Ex. 4.8(5)]{bib:KL15a}. This property occurs even if
ultracontractivity of $\semi {\tilde T}$ fails to hold, which makes a striking difference with the classical case. As it will be seen below, essential
differences between classical and non-local Schr\"odinger operators also occur in respect of their hyper- and supercontractivity properties.

The main results of the present paper are as follows.
We first work in a general setting with $L^2$-semigroups associated with self-adjoint operators $H$ bounded from below, which satisfy some mild
regularity assumptions (Section 2). In this framework we introduce the concept of $L^p$-\emph{ground state domination} of the semigroup $\{T_t: t\geq 0\}$ and
its asymptotic version (Definition \ref{def:pgsd}). By using them, we refine 
the notions of hypercontractivity and supercontractivity of
the transformed semigroup $\{\widetilde T_t: t\geq 0\}$ (Theorems \ref{thm:ihc_gsd}-\ref{thm:isc_gsd}), and derive necessary and sufficient conditions for
these properties to hold (Proposition \ref{prop:suff}, Theorem \ref{thm:nec}, Corollary \ref{cor:nec_suf} and Theorem \ref{thm:border}). These conditions
are expressed explicitly in terms of the behaviour of the L\'evy density $\nu$ and the potential $V$ at infinity. Surprisingly, they are in match with the conditions
under which ultracontractivity and its asymptotic version hold, recently obtained in \cite{bib:KL15a}, and leads to unexpected equivalences in the smoothing
properties of semigroups of non-local Schr\"odinger operators. For the class of jump-paring L\'evy processes and for a large class of confining potentials,
we obtain that supercontractivity and hypercontractivity of the transformed semigroup $\{\widetilde T_t: t\geq 0\}$ coincide respectively with ultracontractivity
and asymptotic ultracontractivity of the semigroup (Corollaries \ref{cor:isp_iuc}-\ref{cor:ihp_iuc}). This has not been observed before, and it makes a stark
contrast to classical Schr\"odinger operators featuring the Laplacian \cite[p.396]{bib:DS}.

\subsection*{Notation}
We will use the notation $C_i(a,b,...)$, $i = 1,2,...$, for positive constants dependent on parameters $a,b,...$ throughout this paper, while dependence
on the process $\pro X$ and on the dimension $d$ is assumed without being explicitly indicated. Auxiliary constants in the proofs are denoted by $c_i$. We
will also use the notation $f \asymp g$ meaning that there exists a constant $C$ such that $C^{-1} g \leq f \leq Cg$ holds. We denote by $\1(x)$ the indicator
function of $\R^d$.

\section{$L^p$-ground state domination properties}

Throughout this section we assume that $H$ is a self-adjoint operator on $L^2(\R^d,dx)$ and bounded from below, such that
\begin{itemize}
\item[(a)]
$\lambda_0:=\inf\Spec H$ is a non-degenerate eigenvalue, i.e., a \emph{ground state} of $H$ exists
\item[(b)]
$e^{-tH}$ is positivity improving for every $t>0$, i.e., $e^{-tH} f(x) > 0$ for all $x \in \R^d$ and $f \in L^2(\R^d,dx)$
such that $f \geq 0$ and $f \neq 0$ a.e.
\item[(c)]
$e^{-tH}$ is bounded on $L^{\infty}(\R^d,dx)$ for all $t > 0$, and satisfies
$$
|e^{-tH} f(x)| \leq e^{-tH} \1(x) \left\|f\right\|_{\infty}, \quad \text{for a.\,e.} \ x \in \R^d \ \text{and every}
\ f \in L^{\infty}(\R^d,dx)
$$
\item[(d)]
there exists $t_{\rm b} >0$ such that $e^{-t_{\rm b} H}$ is bounded from $L^2(\R^d,dx)$ to $L^{\infty}(\R^d,dx)$.
\end{itemize}
The corresponding unique $L^2$-normalized eigenfunction $\varphi_0$, called \emph{ground state}, can be assumed to be strictly
positive \cite[Th. XIII.43]{bib:RS}, \cite[Sect. 3.4.3]{bib:LHB}. Moreover, it immediately follows from the eigenvalue equation
$e^{-t_{\rm b} H} \varphi_0 = e^{-\lambda_0 t_{\rm b}} \varphi_0$ and (d) that $\varphi_0 \in L^{\infty}(\R^d,dx)$.

The following properties of the semigroup $\{e^{-tH}: t \geq 0\}$ will be central to our investigations. Denote by $L^p(\R^d,\mu)$,
$p \in [1,\infty)$, the space of $L^p$-integrable functions on $\R^d$ weighted by the probability measure $\mu(dx)=\varphi_0^2(x)dx$.
Clearly, we have the identification $L^{\infty}(\R^d,\mu) = L^{\infty}(\R^d,dx)$, and for every $1 \leq q <p \leq \infty$ the inclusions
$L^p(\R^d,\mu) \subset L^q(\R^d,\mu)$ hold.

\begin{definition}[\textbf{Ground state domination properties}]
\label{def:pgsd}
Let $p \in (2,\infty]$. We say that
\begin{itemize}
\item[(i)]
the operator $e^{-tH}$ is \emph{$L^p$-ground state dominated} (abbreviated as $L^p$-GSD) if
\begin{align} \label{tpgsd}
\frac{e^{-tH} \1}{\varphi_0} \in L^p(\R^d, \mu)
\end{align}
\item[(ii)]
the semigroup $\{e^{-tH}: t \geq 0\}$ is \emph{$L^p$-ground state dominated} (abbreviated as $L^p$-GSD) if for every
$t>0$ the operators $e^{-tH}$ are $L^p$-ground state dominated
\item[(iii)]
the semigroup $\{e^{-tH}: t \geq 0\}$ is \emph{asymptotically $L^p$-ground state dominated} (abbreviated as $L^p$-AGSD)
if there exists $t_p >0$ such that for every $t>t_p$ the operators $e^{-tH}$ are $L^p$-ground state dominated. If the specific
value of $t_p$ is essential, we write $(t_p,L^p)$-AGSD to emphasize this.
\end{itemize}
\end{definition}
\noindent
For $p =\infty$ such a domination property has been recently introduced and used in \cite{bib:KL15a}. The definition above
considers this now for all $p \in (2,\infty]$.
The $L^p$-GSD property expresses (in terms of appropriate weighted $L^p$-spaces) the balance between the fall-off of the mass
of the semigroup and the decay of the ground state at infinity. The ground state $\varphi_0$ is a key concept in many applications
in quantum and statistical physics (see, e.g., \cite{bib:KL15d}), thus obtaining further ``multi-scale" information according to
the $L^p$-norms on its regularity and localization properties, which is our goal in this paper, is clearly of interest. 

We define the \emph{ground state-transformed semigroup} $\{e^{-t\widetilde{H}}: t \geq 0\}$ by
$$
e^{-t\widetilde{H}} f(x) = \frac{e^{\lambda_0 t}}{\varphi_0(x)} e^{-tH}(f \varphi_0)(x), \quad f \in L^2(\R^d, \mu), \ \ \ t \geq 0.
$$
The following observations will be fundamental in what follows.
\begin{lemma}
\label{lem:pgsd}
Let $p \in (2,\infty)$. Consider the following two conditions.
\begin{itemize}
\item[(1)] For some $t >0$ the operator $e^{-tH}$ is $L^p$-GSD.
\item[(2)] For some $t > 0$ the operator $e^{-t \widetilde{H}}$ is bounded from $L^2(\R^d,\mu)$ to $L^p(\R^d,\mu)$.
\end{itemize}
Then we have the following:
\begin{itemize}
\item[(i)]
If (1) holds for some $t=s > 0$, then (2) follows for $t = s + t_{\rm b}$, where $t_{\rm b}$ is given by (d).
\item[(ii)]
If (2) holds for some $t = s >0$ and
\begin{align} \label{eq:gs_int}
\varphi_0^{1-\frac{1}{p-1}} \in L^1(\R^d,dx),
\end{align}
then (1) follows for $t = 2s$.
\end{itemize}
\end{lemma}

\begin{proof}
We first prove (i). We have for all $f \in L^2(\R^d,\mu)$,
\begin{align*}
\int_{\R^d} | e^{-(t+t_{\rm b})\widetilde{H}} f(x)|^p \mu(dx) = \int_{\R^d} \left(\frac{e^{\lambda_0(t+t_{\rm b})}}{\varphi_0(x)}
\left|e^{-tH} e^{-t_{\rm b}H} (f \varphi_0)(x)\right|\right)^p \mu(dx).
\end{align*}
By the standing assumptions (c) and (d), we have
$$
\left|e^{-tH} e^{-t_{\rm b}H} (f \varphi_0)(x)\right| \leq e^{-tH} \1(x) \left\|e^{-t_{\rm b}H} (f \varphi_0)\right\|_{\infty} \leq
e^{-tH} \1(x) \left\|e^{-t_{\rm b}H}\right\|_{2,\infty} \left\|f \right\|_{L^2(\R^d,\mu)},
$$
for almost all $x \in \R^d$. Thus
$$
\left(\int_{\R^d} |e^{-(t+t_{\rm b})\widetilde{H}} f(x)|^p \mu(dx)\right)^{1/p} \leq
\left\|e^{-t_{\rm b}H}\right\|_{2,\infty} e^{\lambda_0(t+t_{\rm b})} \left(\int_{\R^d} \left|\frac{ e^{-tH} \1(x)}{\varphi_0(x)}\right|^p
\mu(dx)\right)^{1/p}
\, \left\|f \right\|_{L^2(\R^d,\mu)},
$$
which completes the proof of (i).

To show (ii), choose $q$ such that $1/p + 1/q = 1$. By the symmetry of the operator $e^{-t\widetilde{H}}$ and the duality of $L^p(\R^d,\mu)$
and $L^q(\R^d,\mu)$, boundedness of $e^{-t\widetilde{H}}$ from $L^2(\R^d,\mu)$ to $L^p(\R^d,\mu)$ implies its boundedness $L^q(\R^d,\mu)$ to
$L^2(\R^d,\mu)$ with the same norm. By \eqref{eq:gs_int} we have $f:=\frac{1}{\varphi_0} \in L^q(\R^d,\mu)$, and thus $g_t :=
e^{-t\widetilde{H}}f \in L^2(\R^d,\mu)$. With this we have
\begin{align*}
\left(\int_{\R^d} \left|\frac{e^{-2tH} \1 (x)}{\varphi_0(x)}\right|^p \mu(dx)\right)^{1/p}
& =
e^{-2\lambda_0 t} \left(\int_{\R^d} \left|e^{-t\widetilde{H}} g_t (x)\right|^p \mu(dx)\right)^{1/p}
 \leq
C_{2,p,t} e^{-2\lambda_0 t}  \left\|g_t\right\|_{L^2(\R^d,\mu)} \\
& \leq C^2_{2,p,t} e^{-2\lambda_0 t}  \left\|f \right\|_{L^q(\R^d,\mu)}
= C^2_{2,p,t} e^{-2\lambda_0 t} \int_{\R^d} \varphi_0^{1-\frac{1}{p-1}}(x) dx < \infty,
\end{align*}
where $C_{2,p,t}$ is the norm of the operator $e^{-t\widetilde{H}}$ from $L^2(\R^d,\mu)$ to $L^p(\R^d,\mu)$.
\end{proof}

\begin{remark} \label{rem:pgsd} {\rm

\noindent

\begin{itemize}
\item[(1)]
An implication of (ii) in Lemma \ref{lem:pgsd} is that whenever
\begin{align} \label{eq:gs_int_delta}
\text{there exists \ $\delta \in (0,1)$ \ such that $\varphi_0^{1-\delta} \in L^1(\R^d,dx)$,}
\end{align}
and $p \geq 1+ 1/\delta$, then (2) implies (1) for appropriate $t$.
\item[(2)]
Condition \eqref{eq:gs_int_delta} (or \eqref{eq:gs_int}) appears to be non-standard. However, often (and in most cases of direct interest,
see e.g. \cite{bib:Car, bib:KL15a, bib:CMS}) they hold for both classical and non-local Schr\"odinger operators. In particular,
\eqref{eq:gs_int_delta} holds for all the examples discussed in \cite[Sect. 4]{bib:KL15a}.

\item[(3)]
It is straightforward to check that $L^p$-GSD implies
\begin{align} \label{eq:gs_int_L1}
\varphi_0 \in L^1(\R^d,dx).
\end{align}
In this light, it is reasonable to ask whether the assumption \eqref{eq:gs_int_delta} in Lemma \ref{lem:pgsd} (ii) could be relaxed.
It does not seem to be immediate to settle if the same is true under the weaker condition \eqref{eq:gs_int_L1} in the generality of
the present framework. In Section \ref{sec:IC} below we will show that for a wide selection of non-local Schr\"odinger operators with
confining potentials for which \eqref{eq:gs_int_L1} automatically holds, the restriction \eqref{eq:gs_int_delta} can be lifted. This
improvement will be based on the observation that assuming \eqref{eq:gs_int_L1} instead of \eqref{eq:gs_int} and by following
through the argument in the proof of Lemma \ref{lem:pgsd} (ii) with $f:=1/\varphi_0^{1/q}$, we have
$$
\frac{e^{-2t H} \varphi_0^{1/p}}{\varphi_0} \in L^p(\R^d,\mu).
$$
This property is much weaker than \eqref{tpgsd}, but still can be applied well whenever strong enough estimates of $\varphi_0$ at
infinity are available.
\end{itemize}
}
\end{remark}

In what follows, the above $L^p$-ground state domination properties will be applied to characterize the hypercontractivity and
supercontractivity properties of transformed semigroups corresponding to non-local Schr\"odinger operators. Recall the following
terminology from \cite{bib:DS}.

\begin{definition}[\textbf{Contractivity properties}]

\noindent
\begin{itemize}
\item[(i)]
The semigroup $\{e^{-t\widetilde{H}}: t \geq 0\}$ is \emph{supercontractive} if for every $p \in (2,\infty)$ and $t   > 0$ the
operators $e^{-t\widetilde{H}}$ are bounded from $L^2(\R^d,\mu)$ to $L^p(\R^d,\mu)$.
\item[(ii)]
The semigroup $\{e^{-tH}: t \geq 0\}$ is called \emph{intrinsically supercontractive} (abbreviated as ISC) if the semigroup
$\{e^{-t\widetilde{H}}: t \geq 0\}$ is supercontractive.
\item[(iii)]
The semigroup $\{e^{-t\widetilde{H}}: t \geq 0\}$ is \emph{hypercontractive} if for every $p \in (2,\infty)$ there exists
$t_p > 0$ such that for every $t \geq t_p$ the operators $e^{-t\widetilde{H}}$ are bounded from $L^2(\R^d,\mu)$ to $L^p(\R^d,\mu)$.
\item[(iv)]
The semigroup $\{e^{-tH}: t \geq 0\}$ is called \emph{intrinsically hypercontractive} (abbreviated as IHC) if the semigroup
$\{e^{-t\widetilde{H}}: t \geq 0\}$ is hypercontractive.
\end{itemize}
\end{definition}

The next two theorems are the main results of this section, providing a necessary and sufficient condition for ISC and IHC in terms
of the $L^p$-GSD and $L^p$-AGSD properties. This will allow us to refine, study and interpret strong smoothness properties of the semigroups
$\{e^{-t\widetilde{H}}: t \geq 0\}$ through the fall-off properties of the mass of $e^{-tH}$, i.e., the decay of the functions
$e^{-tH} \1(x)$ when $|x| \to \infty$.

\begin{theorem}  \emph{\textbf{(IHC and $L^p$-AGSD)}}
\label{thm:ihc_gsd}
The following hold:
\begin{itemize}
\item[(i)]
If the semigroup $\{e^{-tH}: t \geq 0\}$  is IHC and there exists $\delta \in (0,1)$ such that $\varphi_0^{1-\delta} \in L^1(\R^d,dx)$,
then for every $p \in (2,\infty)$ the semigroup $\{e^{-tH}: t \geq 0\}$ is $L^p$-AGSD.
Specifically, if for every $p \in (2,\infty)$ there exists $t_p>0$ such that for all
$t \geq t_p$ the operators $e^{-t\widetilde{H}}$ are bounded from $L^2(\R^d,\mu)$ to $L^p(\R^d,\mu)$,
then for every $p \in (2,\infty)$ the semigroup $\{e^{-tH}: t \geq 0\}$ is $\big(2(t_p \vee t_{p_{\delta}}),L^p\big)$-AGSD,
where $p_{\delta} = 1+ 1/\delta$.

\item[(ii)]
If for every $p \in (2,\infty)$ the semigroup $\{e^{-tH}: t \geq 0\}$ is $L^p$-AGSD, then the semigroup $\{e^{-tH}: t \geq 0\}$ is IHC.
Specifically, if for every $p \in (2,\infty)$ there exists $t_p>0$ such that for all $t \geq t_p$ the semigroup  $\{e^{-tH}: t \geq 0
\}$ is $(t_p,L^p)$-AGSD, then for every $p \in (2,\infty)$ and all $t \geq t_{\rm b}+t_p$ the operators $e^{-t\widetilde{H}}$ are bounded
from $L^2(\R^d,\mu)$ to $L^p(\R^d,\mu)$.
\end{itemize}
\end{theorem}

\begin{proof}
This is obtained by a direct application of Lemma \ref{lem:pgsd} (and Remark \ref{rem:pgsd} (1)) and the monotonicity in $p \in (2,\infty)$
of the inclusions of the $L^p(\R^d,\mu)$ spaces.
\end{proof}

\begin{theorem}  \emph{\textbf{(ISC and $L^p$-GSD)}}
\label{thm:isc_gsd}
The following hold:
\begin{itemize}
\item[(i)]
If the semigroup $\{e^{-tH}: t \geq 0\}$ is ISC and there exists $\delta \in (0,1)$ such that $\varphi_0^{1-\delta} \in L^1(\R^d,dx)$,
then for every $p \in (2,\infty)$ the semigroup $\{e^{-tH}: t \geq 0\}$ is $L^p$-GSD.
\item[(ii)]
If for every $p \in (2,\infty)$ the semigroup $\{e^{-tH}: t \geq 0\}$ is $L^p$-GSD and for every $t>0$ the operators
$e^{-tH}:L^2(\R^d,dx) \to L^{\infty}(\R^d)$ are bounded (i.e., the semigroup $\big\{e^{-tH}: t \geq 0 \big\}$ is ultracontractive),
then the semigroup $\{e^{-tH}: t \geq 0\}$ is ISC.
\end{itemize}
\end{theorem}

\begin{proof}
We again apply Lemma \ref{lem:pgsd} for all $t>0$.
\end{proof}

\section{Intrinsic contractivity-type properties for jump-paring L\'evy processes} \label{sec:IC}
\subsection{Underlying L\'evy processes and the corresponding non-local Schr\"odinger operators}
\noindent
Let $\pro X$ be a symmetric L\'evy process with values in $\R^d$, $d \geq 1$, with probability measure $\pr^x$
of the process starting from $x \in \R^d$. We use the notation $\ex^x$ for expectation with respect to $\pr^x$.
Recall that $\pro X$ is a Markov process with respect to its natural filtration, satisfying the strong Markov
property and having c\`adl\`ag paths. It is determined by the characteristic function
$$
\ex^0 \left[e^{i \xi \cdot X_t}\right] = e^{-t \psi(\xi)}, \quad \xi \in \R^d, \ t >0,
$$
with the characteristic exponent given by the L\'evy-Khintchin formula
\begin{align} \label{eq:Lchexp}
\psi(\xi) = A \xi \cdot \xi + \int_{\R^d} (1-\cos(\xi \cdot z)) \nu(dz).
\end{align}
Here $A=(a_{ij})_{1\leq i,j \leq d}$ is a symmetric non-negative definite matrix, and $\nu$ is a symmetric
L\'evy measure on $\R^d \backslash \left\{0\right\}$, i.e., $\int_{\R^d} (1 \wedge |z|^2) \nu(dz) < \infty$ and
$\nu(E)= \nu(-E)$, for every Borel set $E \subset \R^d \backslash \left\{0\right\}$. For more details
on L\'evy processes we refer to \cite{bib:Sat,bib:J}.

We will assume throughout that
\begin{align} \label{eq:nuinf}
\nu(\R^d \backslash \left\{0\right\})=\infty \quad \text{and} \quad \nu(dx)=\nu(x)dx, \quad \text{with}
\quad \nu(x)     > 0.
\end{align}
For simplicity, we denote the density of the L\'evy measure also by $\nu$ as it is the object we will use below.
Note that the properties (\ref{eq:nuinf}) jointly imply that $\pro X$ is a strong Feller
process, or equivalently,
there exist measurable transition probability densities $p(t,x,y) = p(t,0,y-x)=: p(t,y-x)$ with respect to Lebesgue measure
such that $\pr^0(X_t \in E) = \int_E p(t,x)dx$, for every Borel set $E \subset \R^d$ (see e.g. \cite[Th. 27.7]{bib:Sat}).
The transition probability densities $p_D(t,x,y)$ of the process killed upon exiting an open bounded set $D \subset \R^d$
are given by the Dynkin-Hunt formula
$$
p_D(t,x,y)= p(t,y-x) - \ex^x\left[ \tau_D < t; p(t-\tau_D, y-X_{\tau_D})\right], \quad x , y \in D,
$$
where $\tau_D = \inf\left\{t \geq 0: X_t \notin D\right\}$ is the first exit time of the process from $D$. The Green function
is given by $G_D(x,y)= \int_0^{\infty} p_D(t,x,y) dt$, for all $x, y \in D$, and $G_D(x,y) = 0$ if $x \notin D$ or $y \notin D$.

We will use the following class of L\'evy processes (cf. \cite{bib:KL15b}).

\begin{definition}[\textbf{Symmetric jump-paring L\'evy processes}]
\label{jumpparing}
Let $\pro X$ be a L\'evy process with L\'evy-Khintchin exponent $\psi$ as in \eqref{eq:Lchexp}--\eqref{eq:nuinf}, satisfying the
following conditions.
\begin{itemize}
\item[\textbf{(A1)}] \emph{\textbf{L\'evy density:}}
There exist a non-increasing profile function $g:(0,\infty) \to (0,\infty)$ and constants $C_1, C_2   >0$ such that
\begin{align} \label{eq:profile_g}
\nu(x) \asymp C_1 g(|x|), \quad x \in \R^d \backslash \left\{0\right\},
\end{align}
and
\begin{align} \label{eq:jp-prop}
\int_{|x-y|>1 \atop |y| > 1} g(|x-y|) g(|y|) dy \leq C_2 \: g(|x|), \quad |x|\geq 1.
\end{align}
\medskip
\item[\textbf{(A2)}] \emph{\textbf{Transition density:}}
There exists $t_{\rm b} >0$ such that $\sup_{x\in \R^d} p(t_{\rm b},x) = p(t_{\rm b},0) < \infty$.
\medskip
\item[\textbf{(A3)}] \emph{\textbf{Green function:}}
For all $0<p<q<R \leq 1$ we have
$$
\sup_{x \in B(0,p)} \sup_{y \in B(0,q)^c} G_{B(0,R)}(x,y) < \infty.
$$
\end{itemize}
We call $\pro X$ satisfying the above conditions a \emph{symmetric jump-paring L\'evy process} and refer to the convolution
condition in \eqref{eq:jp-prop} as the \emph{jump-paring property}.
\end{definition}
\noindent
The bound in \eqref{eq:jp-prop} provides a control of the convolutions
of $\nu$ with respect to large jumps and has a structural importance in defining the class of processes we consider.
It says that the intensity of double large jumps of the process are dominated by the intensity of a single large jump.
Let $\nu_{1}(x) = \nu(x) \1_{B(0,1)^c}(x)$. It is then seen iteratively that under \eqref{eq:jp-prop} in fact
$$
\nu_{1}^{n*}(x) \leq C_3^{n-1} \nu_{1}(x), \quad |x| \geq 1, \; n \in \N,
$$
holds, which means that every sequence of any finite length of large jumps of the process is dominated by single large
jumps, which gives the name to the class of L\'evy processes above.

The convolution condition \eqref{eq:jp-prop} has been introduced in \cite{bib:KL15a} and proved to be a strong tool in studying
large-scale properties of non-local Schr\"odinger semigroups related to jump L\'evy processes. Results obtained in \cite{bib:KL15a}
include sharp estimates on the ground state and upper estimates on the other eigenfunctions at infinity, and necessary and sufficient
conditions for intrinsic ultracontractivity and its asymptotic version obtained via $L^\infty$-GSD and $L^\infty$-AGSD properties. Our
present investigations complement this by focusing on further contractivity properties: $L^p$-GSD for $p \in (2,\infty)$, intrinsic
supercontractivity and intrinsic hypercontractivity. Since below we often refer to \cite{bib:KL15a}, we note that under
\eqref{eq:Lchexp}--\eqref{eq:nuinf}, our present conditions (A2)-(A3) coincide with Assumptions 2.2-2.3 there, while (A1) is a variant
of Assumption 2.1. For simplicity, in the present paper this assumption is stated in terms of a profile function $g$, which gives some
more regularity on the behaviour of $\nu$ around the origin; this is only a slight technical difference which has no impact on the
results obtained here or in the referred paper.
Recently, in \cite{bib:KS14} condition \eqref{eq:jp-prop} has also been used to characterize the short-time behaviour of heat kernels
for a large class of convolution semigroups. It can easily be checked that (A1) in fact implies the comparability condition
\cite[Lem. 1(a)]{bib:KS14}
\begin{align} \label{eq:comparability}
g(|x|) \asymp g(|x|+1), \quad |x| \geq 1.
\end{align}

Assumption (A2) is equivalent with $e^{-t_{\rm b} \psi} \in L^1(\R^d)$, for some $t_{\rm b}   >0$. In this case $p(t_{\rm b},x)$ can be
obtained by the Fourier inversion formula, which extends to all $t \geq t_{\rm b}$ by the Markov property of $\pro X$. Further
details on the existence and properties of transition probability densities for L\'evy processes can be found in \cite{bib:KSch} and
the references therein. We also note that in many cases of interest Assumption (A3) follows directly from space-time estimates of the
densities $p(t,x)$. Indeed, if Assumption (A2) holds and for every $r>0$ there exists $C=C(r)$ such that $\sup_{|x|\geq r} p(t,x) \leq
C t$, $t>0$, then (A3) results by standard estimates, see \cite[Lem. 2.2]{bib:KL15b} and \cite[Prop. 2.3]{bib:BKK}.

The generator $L$ of the process $\pro X$ is uniquely determined by its Fourier symbol
\begin{align} \label{def:gen}
\widehat{L f}(\xi) = - \psi(\xi) \widehat{f}(\xi), \quad \xi \in \R^d, \; f \in D(L),
\end{align}
with domain $D(L)=\{f \in L^2(\R^d): \psi \widehat f \in L^2(\R^d)\}$. The generator is a negative, non-local self-adjoint operator with
core $C_0^{\infty}(\R^d)$, and
$$
L f(x) = \sum_{i,j=1}^d a_{ij} \frac{\partial^2 f}{\partial x_j \partial x_i} (x)
+  \int_{\R^d} \big(f(x+y)-f(x) - y \cdot \nabla f(y) \1_{\left\{|y| \leq 1 \right\}}\big)\nu(y)dy, \quad x \in \R^d,
\; f \in C_0^{\infty}(\R^d).
$$
There is a vast supply of L\'evy processes satisfying Assumptions (A1)-(A3), including large subclasses of isotropic unimodal L\'evy
processes, subordinate Brownian motions, L\'evy processes with non-degenerate Brownian components, symmetric stable-like processes, or
processes with sub-exponentially or exponentially localized L\'evy measures. In particular, it covers \emph{non-Gaussian isotropic stable} and \emph{relativistic stable processes} ($L= - (-\Delta + m^{2/\alpha})^{\alpha/2} + m$, $\alpha  \in (0,2)$, $m \geq 0$), \emph{jump-diffusions} ($L= \Delta - (-\Delta)^{\alpha/2}$, $\alpha  \in (0,2)$), \emph{geometric stable processes} ($L = -\log(1+(-\Delta)^{\alpha/2})$, $\alpha \in (0,2)$) and all of the examples discussed in detail in \cite[Sect. 4]{bib:KL15a}.

We choose the class of potentials in a way which allows us to construct Feynman-Kac semigroups.
\begin{definition}[\textbf{$X$-Kato class}]
{\rm
We say that the Borel function $V: \R^d \to \R$ called \emph{potential} belongs to \emph{Kato-class} $\cK^X$ associated
with the L\'evy process $\pro X$ if it satisfies
\begin{align}
\label{eq:Katoclass}
\lim_{t \downarrow 0} \sup_{x \in \R^d} \ex^x \left[\int_0^t |V(X_s)| ds\right] = 0.
\end{align}
Also, we say that $V$ is an \emph{$X$-Kato decomposable potential}, 
whenever
$$
V=V_+-V_-, \quad \text{with} \quad V_- \in \cK^X \quad \text{and} \quad V_+ \in \cK^X_{\loc},
$$
where $V_+$, $V_-$ denote the positive and negative parts of $V$, respectively, and where $V_+ \in \cK^X_{\loc}$
means that $V_+ 1_B \in \cK^X$ for all compact sets $B \subset \R^d$.
\label{Xkato}
}
\end{definition}
\noindent
For simplicity, we refer to $X$-Kato decomposable potentials as \emph{$X$-Kato class potentials}. It is straightforward that
$L^{\infty}_{\loc}(\Rd) \subset \cK_{\loc}^X$, and by stochastic continuity of $\pro X$ also $\cK_{\loc}^X \subset L^1_{\loc}(\R^d)$.
Note that condition \eqref{eq:Katoclass} allows local singularities of $V$. For specific processes $\pro X$ the definition of $X$-Kato
class can be reformulated more explicitly in terms of the kernel $p(t,x)$ restricted to small $t$ and small $x$. It is shown in
\cite[Cor. 1.3]{bib:GrzSz} that \eqref{eq:Katoclass} is equivalent with
\begin{align} \label{eq:Kato_new}
\lim_{t \to 0^{+}} \sup_{x \in \R^d} \int_0^t \int_{B(x,t)} p(s,x-y)|V(y)| dy ds = 0.
\end{align}
In this section we consider \emph{confining potentials} in the following sense.

\medskip

\begin{itemize}
\item[\textbf{(A4)}]
Let $V$ be an $X$-Kato class potential such that $V(x) \to \infty$ as $|x| \to \infty$.
\end{itemize}

\medskip

Next we define
$$
T_t f(x) = \ex^x\left[e^{-\int_0^t V(X_s) ds} f(X_t)\right], \quad f \in L^2(\R^d), \ t>0.
$$
Standard arguments based on Khasminskii's Lemma (see, e.g., \cite[Lem. 3.37-3.38]{bib:LHB}) imply for an $X$-Kato class potential
$V$ that there exist constants $C_4, C_5 > 0$ such that
\begin{align}
\label{eq:khas}
\sup_{x \in \R^d} \ex^x\left[e^{-\int_0^t V(X_s)ds}\right] \leq \sup_{x \in \R^d}
\ex^x\left[e^{\int_0^t V_-(X_s)ds}\right] \leq C_4 e^{C_5t}, \quad t>0.
\end{align}
Using the Markov property and stochastic continuity of the process it can be shown that $\{T_t: t\geq 0\}$ is a strongly
continuous semigroup of symmetric operators on $L^2(\R^d)$, which we call the \emph{Feynman-Kac semigroup} associated with the
process $\pro X$ and potential $V$. In particular, by the Hille-Yoshida theorem there exists a self-adjoint operator $H$, bounded
from below, such that $e^{-t H} = T_t$. We call the operator $H$ a \emph{non-local Schr\"odinger operator} based on the infinitesimal
generator $L$ of the process $\pro X$. Since any $X$-Kato class potential is relatively form bounded with respect to $H_0=-L$ with
relative bound less than 1, the operator $H=H_0+V$ is also well-defined in form sense, see \cite[Ch. 2]{DC}.

The following are some basic properties of the operators $T_t$ needed below, for a proof see \cite[Lem. 2.1]{bib:KL15a}.

\begin{lemma}
\label{lm:semprop}
Let $\pro X$ be a symmetric L\'evy process with L\'evy-Khintchin exponent satisfying \eqref{eq:Lchexp}-\eqref{eq:nuinf}
such that Assumption (A2) holds, and let $V$ be an $X$-Kato class potential. Then the following properties hold:
\begin{itemize}
\item[(1)]
For all $t>0$, every $T_t$ is a bounded operator on every $L^p(\R^d,dx)$ space, $1 \leq p \leq \infty$. The operators
$T_t: L^p(\R^d,dx) \to L^p(\R^d,dx)$ for $1 \leq p \leq \infty$, $t > 0$, and $T_t: L^p(\R^d,dx) \to L^{\infty}(\R^d,dx)$
for $1 < p \leq \infty$, $t \geq t_{\rm b}$, and $T_t: L^1(\R^d,dx) \to L^{\infty}(\R^d,dx)$ for $t \geq 2t_{\rm b}$
are bounded.
\vspace{0.1cm}
\item[(2)]
For all $t  >0$ the operators $T_t:L^2(\R^d,dx) \to L^2(\R^d,dx)$ are compact.
\item[(3)]
For all $t>0$ and $f \in L^{\infty}(\R^d,dx)$, $T_t f$ is a bounded continuous function.
\vspace{0.1cm}
\item[(4)]
For all $t > 0$ the operators $T_t$ are positivity improving.
\end{itemize}
\end{lemma}
\noindent
Note that we do not
assume that $p(t,x)$ is bounded for all $t>0$, and thus in general the operators $T_t: L^2(\R^d,dx) \to L^{\infty}(\R^d,dx)$
need not be bounded for $t<t_{\rm b}$. Also, note that from the above it follows that the semigroup $\{T_t: t\geq 0\}$ satisfies
the basic regularity conditions (a)-(d) in Section 2.

The following upper bound on the ground state of $H$ at infinity obtained in \cite[Cor. 2.2]{bib:KL15a} will be an important
ingredient below, guaranteeing that for the class of processes and potentials considered in this section we have $\varphi_0
\in L^1(\R^d,dx)$. For $r> 0$ denote $V^{*}_r(x):=\sup_{y \in B(x,r)} V(y)$, $x \in \R^d$.

\begin{proposition} \label{prop:up_gs}
Let $\pro X$ be a symmetric L\'evy process determined by \eqref{eq:Lchexp}-\eqref{eq:nuinf} with defining parameters $A$ and $\nu$
such that Assumptions (A1)-(A3) hold, and let $V$ be a potential satisfying Assumption (A4). Then for every $r \in (0,1]$
there exists $C_6 >0$ and $R>0$ such that
$$
\frac{1}{C_6} \frac{\nu(x)}{V^{*}_r(x)} \leq \varphi_0(x) \leq C_6 \nu(x), \quad |x| \geq R.
$$
\end{proposition}

\smallskip

\subsection{Necessary and sufficient conditions for $L^p$-ground state domination}

We begin with the following result providing sufficient conditions for $L^p$-GSD and $L^p$-AGSD. It directly follows from
the fact that $L^\infty$-GSD and $L^\infty$-AGSD imply $L^p$-GSD and $L^p$-AGSD for any $p \in (2,\infty)$, respectively.

\begin{proposition} \emph{\textbf{(Sufficient conditions for $L^p$-GSD and $L^p$-AGSD)}} \label{prop:suff}
Let $\pro X$ be a symmetric L\'evy process determined by \eqref{eq:Lchexp}-\eqref{eq:nuinf} with defining parameters $A$ and $\nu$ such
that Assumptions (A1)-(A3) hold, and let $V$ be a potential satisfying Assumption (A4). Then the following hold:
\begin{itemize}
\item[(i)]
If there exist constants $C_7>0$ and $R> 0$ such that
$$
\frac{V(x)}{|\log \nu(x)|} \geq C_7, \quad \text{for every \ $|x| \geq R$,}
$$
then for every $p \in (2,\infty)$ the semigroup $\{T_t:t \geq 0\}$ is $(t_0,L^p)$-AGSD with $t_0 = 4/C_7$.
\item[(ii)]
If
$$
\lim_{|x| \to \infty} \frac{V(x)}{|\log \nu(x)|} = \infty,
$$
then for every $p \in (2,\infty)$ the semigroup $\{T_t:t \geq 0\}$ is $L^p$-GSD.
\end{itemize}
\end{proposition}

\begin{proof}
By \cite[Ths. 2.6-2.7]{bib:KL15a} the conditions in (i) and (ii) imply $(t_0,L^\infty)$-AGSD (with $t_0=4/C_7$) and $L^\infty$-GSD,
respectively. Moreover, as seen above, $L^{\infty}(\R^d,\mu) \subset L^p(\R^d,\mu)$, for any $p \in (2,\infty)$.
\end{proof}

It is tempting to expect that $L^\infty$-GSD and $L^\infty$-AGSD are much stronger than $L^p$-GSD and $L^p$-AGSD for $p <\infty$, and so
the above proposition may seem not to give a sharp result. However, this intuition is false. Below we prove that for a large class of
confining potentials $L^\infty$-GSD and $L^\infty$-AGSD are equivalent with $L^p$-GSD and $L^p$-AGSD, for every $p \in (2,\infty)$.

The following theorem is our main result in this section, giving direct necessary conditions for $L^p$-GSD and $L^p$-AGSD. This observation
actually leads to a full characterization of these properties in terms of the L\'evy density $\nu$ and the potential $V$.

\begin{theorem} \emph{\textbf{(Necessary conditions for $L^p$-GSD and $L^p$-AGSD)}} \label{thm:nec}
Let $\pro X$ be a symmetric L\'evy process determined by \eqref{eq:Lchexp}-\eqref{eq:nuinf} with defining parameters $A$ and $\nu$ such
that Assumptions (A1)-(A3) hold, and let $V$ be a potential satisfying (A4). Then we have:
\begin{itemize}
\item[(i)]
If for some $p \in (2,\infty)$ the semigroup $\{T_t:t\geq0\}$ is $(t_0,L^p)$-AGSD, then for every $r \in (0,1)$ and $\varepsilon \in
(0,(p-2)/(pt_0))$ there exists $R >0$ such that
$$
\frac{V^{*}_r(x)}{|\log \nu(x)|} \geq \frac{p-2}{p \, t_0} -\varepsilon, \quad \text{for every \ $|x| \geq R$.}
$$
\item[(ii)]
If for some $p \in (2,\infty)$ the semigroup $\{T_t:t\geq0\}$ is $L^p$-GSD, then for every $r \in (0,1)$
$$
\lim_{|x| \to \infty} \frac{V^{*}_r(x)}{|\log \nu(x)|} = \infty.
$$
\end{itemize}
\end{theorem}

\begin{proof}
(i) By the definition of $(t_0,L^p)$-AGSD, for every $x \in \R^d$ and $r \in (0,1]$ we have
\begin{align} \label{eq:eq11}
\int_{|x-y| < r/2} \frac{|T_{t_0} \1(y)|^p}{\varphi_0^{p-2}(y)} dy \leq \int_{\R^d} \left|\frac{T_{t_0} \1(y)}{\varphi_0(y)} \right|^p
\varphi_0^2(y) dy < \infty
\end{align}
We first estimate the term under the integral at the left hand side above for $y \in B(x,r/2)$:
\begin{align*}
T_{t_0} \1(y) &  \geq \ex^y\left[e^{-\int_0^{t_0} V(X_s)ds}; t_0 < \tau_{B(y,r/2)}\right] \\ &
\geq
e^{- t_0 \sup_{z \in B(y,r/2)}V(z)} \pr^y(t_0 < \tau_{B(y,r/2)}) \geq e^{- t_0 V_r^{*}(x)} \pr^0(t_0 < \tau_{B(0,r/2)}).
\end{align*}
By Proposition \ref{prop:up_gs} and \eqref{eq:comparability} there exists $R \geq 1$ such that
$$
\varphi_0(y) \leq c_1 \nu(y)  \leq c_2 \nu(x) \quad \text{and} \quad \nu(x) \leq 1, \quad \text{whenever \ $|x|   > R$ \ and \ $|x-y|
\leq r/2$}.
$$
Thus by \eqref{eq:eq11} there exist constants $c_3=c_3(p,t_0)$ and $c_4=c_4(p,t_0,r)$ such that for $|x| > R$
$$
e^{- t_0 V_r^{*}(x)} \leq c_3 \frac{\nu(x)^{\frac{p-2}{p}}}{(|B(0,r/2)|)^{1/p} \pr^0(t_0 < \tau_{B(0,r/2)})}
$$
and
$$
\frac{V_r^{*}(x)}{|\log \nu(x)|} \geq \frac{p-2}{p \, t_0} - \frac{c_4}{|\log \nu(x)|}.
$$
Since $|\log \nu(x)| \to \infty$ as $|x| \to \infty$, this completes the proof of (i).

To show (ii), it suffices to observe that $L^p$-GSD implies $(t_0,L^p)$-AGSD for any $t_0   >0$. By following through the above argument,
we obtain
$$
\liminf_{|x| \to \infty} \frac{V_r^{*}(x)}{|\log \nu(x)|} \geq \frac{p-2}{p \, t_0}.
$$
Letting $t_0 \to 0^{+}$, the claim follows.
\end{proof}

It follows from Theorems \ref{thm:ihc_gsd}-\ref{thm:isc_gsd} that under condition \eqref{eq:gs_int_delta} Theorem \ref{thm:nec} above
also gives necessary conditions for IHC and ISC. Next we show that in the framework of this section the same holds without requiring
\eqref{eq:gs_int_delta}.

\begin{corollary} \emph{\textbf{(Necessary and sufficient conditions for ISC and IHC)}} \label{cor:nec_suf}
Let $\pro X$ be a symmetric L\'evy process determined by \eqref{eq:Lchexp}-\eqref{eq:nuinf} with defining parameters $A$ and $\nu$
such that Assumptions (A1)-(A3) hold, and let $V$ be a potential satisfying Assumption (A4). Then we have:
\begin{itemize}
\item[(i)]
The condition in Proposition \ref{prop:suff} (i) is sufficient for IHC. Conversely, if IHC holds, then for every $p \in (3,\infty)$ and
$r \in (0,1)$ there exist $t_p >0$ and $R >0$ such that
$$
\frac{V^{*}_r(x)}{|\log \nu(x)|} \geq \frac{p-3}{2 p t_0}, \quad \text{for every \ $|x| \geq R$.}
$$
\item[(ii)]
The condition in Proposition \ref{prop:suff} (ii) is sufficient for ISC. Conversely, if ISC holds, then for every $r \in (0,1)$
$$
\lim_{|x| \to \infty} \frac{V^{*}_r(x)}{|\log \nu(x)|} = \infty.
$$
\end{itemize}
\end{corollary}

\begin{proof}
We only need to consider the necessary conditions. Similarly as before, it is enough to justify that in (i). Suppose IHC holds and
fix $r \in (0,1]$. By \eqref{eq:gs_int_delta} and the Remark \ref{rem:pgsd} (3) we have that for every $p \in (3,\infty)$ there exists
$t_p >0$ such that $\frac{e^{-2t_p H} \varphi_0^{\frac{1}{p}}}{\varphi_0} \in L^p(\R^d,\mu)$. Thus by using both the lower and upper
bounds in Proposition \ref{prop:up_gs} and by following the argument in the proof of Theorem \ref{thm:nec} (i), we get that there
exists $R >0$ such that for $|x| > R$
$$
\frac{e^{- t_p V_r^{*}(x)}}{\big(V_r^{*}(x)\big)^{1/p}} \leq c_1 \frac{\nu(x)^{\frac{p-3}{p}}}{(|B(0,r/2)|)^{1/p} \pr^0(t_p < \tau_{B(0,r/2)})}
$$
and
$$
\frac{V_r^{*}(x) + (1/p)\log V_r^{*}(x)}{|\log \nu(x)|} \geq \frac{p-3}{p \, t_p} - \frac{c_2}{|\log \nu(x)|},
$$
with some constants $c_1, c_2   >0$. We see that by increasing $R >0$ if necessary, we get the claimed inequality.
\end{proof}

Then the following characterization is a direct corollary of Proposition \ref{prop:suff} and Theorem \ref{thm:nec}.
\begin{theorem} \emph{\textbf{(Characterization of $L^p$-GSD and $L^p$-AGSD)}} \label{thm:border}
Let $\pro X$ be a symmetric L\'evy process determined by \eqref{eq:Lchexp}-\eqref{eq:nuinf} with defining parameters $A$ and
$\nu$ such that (A1)-(A3) hold, and let $V$ satisfy Assumption (A4). Moreover, suppose that the potential $V$ satisfies at least
one of the following additional assumptions:
\begin{itemize}
\item[\textbf{(A5)}]
There exist $r \in (0,1)$ and $R > 0$ such that $V_r^{*}(x) \asymp V(x)$, for $|x| > R$.
\item[\textbf{(A6)}]
There exist a non-decreasing function $f$ and $R>0$ such that $V(x) \asymp f(|x|)$, for $|x| > R$.
\end{itemize}
The following hold:
\begin{itemize}
\item[(i)]
The semigroup $\{T_t:t \geq 0\}$ is $L^p$-GSD for every $p \in (2,\infty)$ if and only if
$$
\lim_{|x| \to \infty} \frac{V(x)}{|\log \nu(x)|} = \infty.
$$
\item[(ii)]
The semigroup $\{T_t:t \geq 0\}$ is $L^p$-AGSD for every $p \in (2,\infty)$ if and only if there exist constants $C>0$ and $R> 0$
such that
$$
\frac{V(x)}{|\log \nu(x)|} \geq C, \quad \text{for every \ $|x| \geq R$.}
$$
\end{itemize}
\end{theorem}

\begin{proof}
Under (A5) the result follows directly from Proposition \ref{prop:suff} and Theorem \ref{thm:nec} above. Suppose now that (A6) holds. It
suffices to prove that the conditions on $V$ and $\nu$ in (i) and (ii) are in fact necessary for $L^p$-GSD and $L^p$-AGSD, respectively.

Consider first (ii). If $(t_0,L^p)$-AGSD holds for some $p \in (2,\infty)$ and $t_0 >0$, then by a similar argument as in the proof of
Theorem \ref{thm:nec} (i) and by (A6), there exist constants $R >0$, $c_1=c_1(p,t_0)$ and $c_2 >0$ (independent of $t_0$ and $p$) such
that
$$
c_2 V(x + x/|x|) \geq \frac{p-2}{p \, t_0} |\log \nu(x)| - c_1, \quad \text{for every \ $|x| \geq R$,}
$$
which is equivalent with
$$
V(x) \geq \frac{p-2}{c_2 p \, t_0} |\log \nu(x-x/|x|)| - \frac{c_1}{c_2}, \quad \text{for every \ $|x| \geq R+1$.}
$$
By (A1) and \eqref{eq:comparability} there exists a constant $c_3 > 0$ (also independent of $t_0$ and $p$)  such that
$$
|\log \nu(x-x/|x|)| \geq |\log \nu(x)| - c_3, \quad \text{for sufficiently large $|x|$.}
$$
Since $|\log \nu(x)| \to \infty$ as $|x| \to \infty$, the assertion (ii) follows. To get (i), it is enough to observe that
$$
\liminf_{|x| \to \infty} \frac{V(x)}{|\log \nu(x)|} \geq  \frac{p-2}{c_2 p \, t_0}
$$
and take the limit $t_0 \to 0^{+}$.
\end{proof}

\medskip
\subsection{Equivalence of intrinsic contractivity properties for non-local Schr\"odinger operators}

The next two results summarize our investigations in the present paper. They say that for a large class of potentials (at least those
satisfying Assumption (A5) or (A6) above) ISC is equivalent to intrinsic ultracontractivity (IUC), and IHS is equivalent to asymptotic
intrinsic ultracontractivity (AIUC). This surprising fact has not been noted before, and it shows a different behaviour from the classical
case. Indeed, for Schr\"odinger operators $H = -\Delta + V$ all these contractivity properties are different, even when the potentials are
quite regular (see \cite[p.336]{bib:DS} and the references therein). First we recall the following definitions from \cite{bib:DS} and
\cite{bib:KL15a}. In what follows, $\{\widetilde T_t: t \geq 0\}$ will denote the ground state-transformed semigroup to $\{T_t: t \geq 0\}$.

\begin{definition} \emph{\textbf{(Intrinsic ultracontractivity-type properties)}}

\noindent
\begin{itemize}
\item[(i)]
The semigroup $\{\widetilde T_t: t \geq 0\}$ is \emph{ultracontractive} if for every $t > 0$ the operators $\widetilde T_t$ are
bounded from $L^2(\R^d,\mu)$ to $L^{\infty}(\R^d)$.
\item[(ii)]
The semigroup $\{T_t: t \geq 0\}$ is called \emph{intrinsically ultracontractive} (abbreviated as IUC) if the semigroup
$\{\widetilde T_t: t \geq 0\}$ is ultracontractive.
\item[(iii)]
The semigroup $\{\widetilde T_t: t \geq 0\}$ is \emph{asymptotically ultracontractive} if there exists $t_0 > 0$ such that for
every $t \geq t_0$ the operators $\widetilde T_t$ are bounded from $L^2(\R^d,\mu)$ to $L^{\infty}(\R^d)$.
\item[(iv)]
The semigroup $\{T_t: t \geq 0\}$ is \emph{asymptotically intrinsically ultracontractive} (abbreviated as AIUC) if the semigroup
$\{\widetilde T_t: t \geq 0\}$ is \emph{asymptotically ultracontractive}.
\end{itemize}
\end{definition}

The following results are direct consequences of the above results and \cite[Th. 2.5, Cor. 2.3]{bib:KL15a}.

\begin{corollary} \emph{\textbf{(Equivalence of intrinsic contractivity properties)}} \label{cor:isp_iuc}
Let $\pro X$ be a symmetric L\'evy process determined by \eqref{eq:Lchexp}-\eqref{eq:nuinf} with defining parameters $A$ and $\nu$ such
that Assumptions (A1)-(A3) hold, and let $V$ be a potential satisfying Assumption (A4) and at least one of the additional Assumptions (A5)
or (A6).  Moreover, suppose that the transition densities $p(t, \cdot)$ are bounded for all $t>0$. Then the following properties are
equivalent:
\begin{itemize}
\item[(i)] $\lim_{|x| \to \infty} \frac{V(x)}{|\log \nu(x)|} = \infty$.
\item[(ii)] The semigroup $\{T_t:t \geq 0\}$ is $L^p$-GSD for every $p \in (2,\infty]$.
\item[(iii)] The semigroup $\{T_t:t \geq 0\}$ is ISC.
\item[(iv)] The semigroup $\{T_t:t \geq 0\}$ is IUC.
\end{itemize}
\end{corollary}

\begin{corollary} \emph{\textbf{(Equivalence of asymptotic intrinsic contractivity properties)}} \label{cor:ihp_iuc}
Let $\pro X$ be a symmetric L\'evy process determined by \eqref{eq:Lchexp}-\eqref{eq:nuinf} with defining parameters $A$ and $\nu$ such that
Assumptions (A1)-(A3) hold, and let $V$ be a potential satisfying Assumption (A4) and at least one of the additional Assumptions (A5) or (A6).
Then the following properties are equivalent:
\begin{itemize}
\item[(i)] There exist  $C, R> 0$ such that $\frac{V(x)}{|\log \nu(x)|} \geq C$, for every $|x| \geq R$.
\item[(ii)] The semigroup $\{T_t:t \geq 0\}$ is $L^p$-AGSD for every $p \in (2,\infty]$.
\item[(iii)] The semigroup $\{T_t:t \geq 0\}$ is IHC.
\item[(iv)] The semigroup $\{T_t:t \geq 0\}$ is AIUC.
\end{itemize}
\end{corollary}

We close this section by noting that for our class of jump-paring L\'evy processes and for confining potentials satisfying (A5) or (A6)
the function $|\log \nu(x)|$ in fact determines the \emph{borderline growth} of the potential for intrinsic contractivity properties listed
in Corollaries \ref{cor:isp_iuc}-\ref{cor:ihp_iuc} above.

\section{Examples}

To illustrate these results, we assume for simplicity that $\pro X$ is a symmetric L\'evy process with characteristic exponent $\psi$
satisfying \eqref{eq:Lchexp}-\eqref{eq:nuinf}, and parameters $A = a \Id$ for $a \geq 0$ and $\nu$ satisfying \eqref{eq:profile_g} for the
profile function $g$, such that $g(r) = r^{-d-\alpha}$ with some $\alpha \in (0,2)$, for all $r \in (0,1]$. Assumptions (A1)-(A3) for
the processes discussed in Examples \ref{eq:eq1} and \ref{eq:eq3} below have been verified in \cite[Sect. 4]{bib:KL15a}; for Example
\ref{eq:eq2} they can be checked similarly, and we leave the details to the reader.

\smallskip

\begin{example}{\rm \textbf{(Processes with polynomially suppressed large jumps)} \label{eq:eq1}
Let $g(r) = r^{-d-\gamma}$ with $\gamma > 0$, for $r \geq 1$. This class includes
\begin{enumerate}
\item[(i)]
\emph{rotationally invariant $\alpha$-stable process} ($\gamma=\alpha$)
\item[(ii)]
\emph{mixture of $i \geq 2$ rotation invariant $\alpha_i$-stable processes with indices $\alpha_i \in (0,2]$} ($\gamma = \min \alpha_i$)
\item[(iii)]
\emph{layered $\alpha$-stable process} ($\gamma > 2$).
\end{enumerate}
Let $V$ be an $X$-Kato class potential such that
$$
V_{+}(x) \asymp f(|x|), \quad \text{with \ $f(r) = (1+r)^{\delta_1} [\log(2+r)]^{\delta_2} [\log(2+\log(2+r)]^{\delta_3}$ \ for
$\delta_1, \delta_2, \delta_3 \in \R$}.
$$
Then the four equivalent conditions in Corollary \ref{cor:isp_iuc} hold if and only if
$$
\text{(1) \ $\delta_1 > 0$, $\delta_2 \in \R$, $\delta_3 \in \R$  \ \ \
or \ \ \ (2) \ $\delta_1 = 0$, $\delta_2> 1$, $\delta_3 \in \R$ \ \ \
or \ \ \ (3) \ $\delta_1 = 0$, $\delta_2=1$, $\delta_3 > 0$.}
$$
Also, the four equivalent conditions in Corollary \ref{cor:ihp_iuc} hold if and only if
$$
\text{(1) \ $\delta_1 > 0$, $\delta_2 \in \R$, $\delta_3 \in \R$  \ \ \
or \ \ \ (2) \ $\delta_1 = 0$, $\delta_2 \geq 1$, $\delta_3 \in \R$ \ \ \
or \ \ \ (3) \ $\delta_1 = 0$, $\delta_2=1$, $\delta_3 \geq 0$.}
$$
}
\end{example}

\smallskip

\begin{example} {\rm \textbf{(Processes with stretched-exponentially suppressed large jumps)} \label{eq:eq2}
Let $g(r) = e^{\frac{c}{\log3}} e^{-c\frac{r^\beta}{\log(2+r)}}$ with $c > 0$, $\beta \in (0,1]$, for $r \geq 1$.
Also, let $V$ be an $X$-Kato class potential such that
$$
V_{+}(x) \asymp f(|x|), \quad \text{with \ $f(r) = (1+r)^{\delta_1} [\log(2+r)]^{\delta_2}$ \ for $\delta_1, \delta_2 \in \R$}.
$$
Then the four equivalent conditions in Corollary \ref{cor:isp_iuc} hold if and only if
$$
\text{(1) \ $\delta_1 > 1$, $\delta_2 \in \R$  \ \ \ \ \ \ \ \
or \ \ \ \ \ \ \ \ (2) \ $\delta_1 = 1$, $\delta_2> -1$.}
$$
Moreover, the four equivalent conditions in Corollary \ref{cor:ihp_iuc} hold if and only if
$$
\text{(1) \ $\delta_1 > 1$, $\delta_2 \in \R$  \ \ \ \ \ \ \ \
or \ \ \ \ \ \ \ \ (2) \ $\delta_1 = 1$, $\delta_2 \geq -1$.}
$$
}
\end{example}

\smallskip

\begin{example} {\rm \textbf{(Processes with exponentially suppressed large jumps)} \label{eq:eq3}
Let $g(r) = e^c e^{-cr }r^{-\gamma}$ with $c > 0$, $\gamma > (d+1)/2$, for $r \geq 1$. This class includes
\begin{enumerate}
\item[(i)]
\emph{relativistic $\alpha$-stable process} ($c=m^{1/\alpha}$, $\gamma=(d+\alpha+1)/2$, $m >0$),
\item[(ii)]
\emph{(exponentially) tempered $\alpha$-stable process} ($c >0$, $\gamma = d+\alpha$).
\end{enumerate}
Let $V$ be an $X$-Kato class potential such that
$$
V_{+}(x) \asymp f(|x|), \quad \text{with \ $f(r) = (1+r)^{\delta_1} [\log(2+r)]^{\delta_2}$ \ for $\delta_1, \delta_2 \in \R$}.
$$
Then the four equivalent conditions in Corollary \ref{cor:isp_iuc} hold if and only if
$$
\text{(1) \ $\delta_1 > 1$, $\delta_2 \in \R$  \ \ \ \ \ \ \ \
or \ \ \ \ \ \ \ \ (2) \ $\delta_1 = 1$, $\delta_2> 0$.}
$$
Moreover, the four equivalent conditions in Corollary \ref{cor:ihp_iuc} hold if and only if
$$
\text{(1) \ $\delta_1 > 1$, $\delta_2 \in \R$  \ \ \ \ \ \ \ \
or \ \ \ \ \ \ \ \ (2) \ $\delta_1 = 1$, $\delta_2 \geq 0$.}
$$
}
\end{example}

\begin{example} {\rm \textbf{(Brownian motion)} \label{eq:eq4}
Consider the Schr\"odinger operator $H =-\Delta +V$ with $V(x) = |x|^{\delta_1} (\log(2+|x|)^{\delta_2}$,
for $\delta_1, \delta_2 \geq 0$, and denote  by $\pro B$ Brownian motion running at twice the usual speed.
As said before, the corresponding Schr\"odinger
semigroup $\{e^{-tH}:t \geq 0\}$ is ISC but not IUC if $\delta_1 = 2$ and $0 < \delta_2 \leq 2$. Note that
$C_{2,\infty,t} = \lim_{p \to \infty} C_{2,p,t} = \infty$, for every $t >0$, where $C_{2,p,t}$ is the norm
of the operator $e^{-t \widetilde{H}}$ from $L^2(\R^d,\mu)$ to $L^p(\R^d,\mu)$, $p \in (2,\infty]$, which shows
a different behaviour from the non-local Schr\"odinger operators studied in this paper
(cf. Corollary \ref{cor:isp_iuc}). Indeed, from the proof
of Corollary \ref{cor:nec_suf} it can be seen that under Assumption (A5) or (A6) the finiteness of $C_{2,p,t}$
for some $p \in (2,\infty)$ in fact guarantees that $C_{2,\infty,t} < \infty$. Moreover, recall that when
$\delta_1 = 2$ and $\delta_2 = 0$, then $\{e^{-tH}:t \geq 0\}$ is IHC but not ISC. This example also shows that for
classical Schr\"odinger operators IHC is not equivalent with AIUC, unlike for non-local Schr\"odinger operators
(cf. Corollary \ref{cor:ihp_iuc}). Using that, see \cite[eq.(2.1)]{MY},
$$
\varphi_0(x) \asymp e^{-\frac{|x|^2}{2}} \quad \text{and} \quad e^{- t H} \1(x) = \ex^x \left[e^{-\int_0^t |B_s|^2 ds} \right]
\asymp e^{-\frac{|x|^2}{2 \coth(2t)}}, \quad x \in \R^d, \ t >0,
$$
we see that for every $p \in (2,\infty)$ there exists $t_p >0$ such that for all $t \geq t_p$ the operator $e^{-tH}$ is $L^p$-GSD.
On the other hand, there is no $t>0$ for which $e^{-tH}$ is $L^{\infty}$-GSD.
}
\end{example}

\end{document}